\documentclass[11pt,twoside,leqno]{aomamlt2e}
\pdfoutput=1

\newcommand{\N}{{\mathbb N}}
\newcommand{\Q}{{\mathbb Q}}
\newcommand{\R}{{\mathbb R}}
\newcommand{\Z}{{\mathbb Z}}
\newcommand{\pz}{{\mathbb P}}

  \pageno{1}

\begin{document}
 \title{Binomial Coefficients and the Distribution of the Primes}
\author{Triantafyllos Xylouris}
Student, University of Bonn/Hannover, September 2007, +49 - (0) 174 - 93 74 007\\
\institution{University of Bonn/Hannover} \email{xtriant@yahoo.com}

\begin{abstract}
Let $\omega (l)=\sum_{p \mid l} 1$ and $m,n \in \N, n \geq m$. We
calculate a formula $\{p \in \pz; p | \binom{l}{k} \}=\pz \cap
\bigcup_i (a(i),b(i)]$ from which we get an identity $\omega
(\binom{nk}{mk}) = \sum_i (\pi (\frac{k}{b(i)}) - \pi
(\frac{k}{a(i)}))+ O(\sqrt{k})$. Erd\ös [Erd79] mentioned that
$\omega (\binom {nk}{mk}) = \log \frac {n^n}{m^m (n-m)^{n-m}}
\frac{k}{\log k} + o(\frac{k}{\log k}).$\\
As an application of the above identities, we conclude some
well-known facts about the distribution of the primes and deduce
$\forall k \in \N$ an expression (also well-known) $\log k = \sum_j
\alpha_k (j)$ which generalizes $\log 2 = \sum_{j=1}^{\infty}
\frac{(-1)^{j+1}}{j}$.
\end{abstract}

\section {Notation}
$\N = \{1,2,3,\dots \} \text{, } \N_0 = \{0,1,2,3, \dots\} \text{, }
 \pz = \{2,3,5,7,11,\dots\} \text{ (primes),}$ $p,p_i,p_{ij}$ will always denote a prime,\\
$\pi (x) = \sum_{p \leq x} 1 \text{ } (x \in \R)$, \\
$e_p(n) = \max \{j \in \N_0 ; p^j \mid n \} \text{ } (p\in \pz , n \in \N) $, \\
$\omega (n) = \sum_{p \mid n} 1 \text{ } (n \in \N)$, we will write $\omega \binom{n}{k}$ for $\omega (\binom{n}{k})$,\\
 $ (m,n)=\max\{d \in \N; d|m \text{ and } d|n\} \text{ } (m,n \in \N)$,\\
 $\psi (x) = \sum_{p^n \leq x} \log p \text{ } (x \in \R)$,\\
$ \log x$ stands for the natural logarithm to the base $e$, we will write $\log ab$ for $\log (ab)$,\\
$x = [x] + \{x\}$, where $[x]=\max\{n \in \Z ; n \leq x\}$ $(x \in \R)$,\\
$(a,b]=\{x \in \R; x \leq b \text{ and } x > a\}$ $(a,b \in \R)$.\\
The $o,O$ - notation will normally be used for $k \longrightarrow \infty$.\\

\section {Introduction}
Obviously, $\binom {2k}{k}$ has every prime number $p \in (k,2k]$ as
a prime divisor. This, together with the fact that $\binom{2k}{k}$
is sufficiently small, is often used for showing $\pi(x) \ll
\frac{x}{\log x}$. On the other hand $\pi(x) \gg \frac{x}{\log x}$
is often shown by noticing that $p^{e_p(\binom{2k}{k})} \leq 2k$ and
$\binom{2k}{k}$ is sufficiently large.

In this paper we somehow sharpen and generalize the above facts for
a wider class of binomial coefficients: Lemma 1 shows exactly which
primes in which intervals divide a given binomial coefficient
$\binom{n}{k}$. Thus, one gets an identity $\omega \binom{n}{k} =
\sum_i (-1)^{a_{n,k}(i)} \pi (\frac{n}{b_{n,k}(i)})$ which we write
down for the binomial coefficients $\binom{nk}{mk}$ ($n,m$ fixed, $k
\rightarrow \infty$) (Theorem 1). \emph{Paul Erd\ös} mentioned an
asymptotic formula for $\omega \binom{n}{k}$ (for $k > n^{1-o(1)}$)
which we write down for $\binom{nk}{mk}$ (Theorem 2, [Erd79]).
Combining Theorem 1 and 2, one gets identities (for all pairs
$(n,m)$) from which one can try to deduce some information about
$\pi(x)$. This is done in Corollary 1 and 2. Corollary 3 is an
application of these identities, which has nothing
to do with primes.\\

Note that the mentioned identities (plus some more) can be derived
significantly easier for $\psi(x)$ in place of $\pi(x)$ (Theorem 3,
e.g. [Land09, pp. 71-95]). The analogues of Corollary 1 and 2 could
then be derived in the same manner for $\psi(x)$. Corollary 3 can
also be proven with those identities which involve $\psi(x)$ instead of $\pi(x)$.\\

\section {The Lemma and the Three Theorems}

Obviously, $\binom {2k}{k}$ has every prime number $p \in (k,2k]$ as
a prime divisor, since these primes divide the numerator $1$ time
and the denominator $0$ times. Lemma 1 is a straightforward
generalization of this fact. For $\binom{2k}{k},
\binom {3k}{k}$ this is partly done in [Fel91].\\
 \begin{description}
 \item[Lemma 1] \emph{Let $n,k \in \N, n \geq k$. We have:}
$$\{p \in \pz; p | \binom{n}{k} \} = \pz \cap( \bigcup_{i=1}^{\infty}
(\left(\bigcup_{j=1}^{\infty} \bigcup_{f=[\frac{n}{k}(j-1)]-j+1
}^{[\frac{n}{k}j]- j - 1} \left(\left(\frac{n-k}{f+1}\right)^{1/i},
\left(\frac{n}{f+j}\right)^{1/i}\right] \right)$$
\begin{equation}\cup \bigcup_{j=1,\text{ } nj \neq 0 \text{ } mod
\text{ } k}^{\infty} \left(
\left(\frac{k}{j}\right)^{1/i},\left(\frac{n}{[\frac{n}{k}j]}\right)^{1/i}\right])).
\end{equation}

\emph{All of the intervals mentioned in (1) are $\neq \emptyset$.
Furthermore, for each fixed $i$, the intervals (depending on $j$ and
$f$, resp. only on $j$) are disjoint.}

 \end{description}
\vspace*{0.5cm} \Proof Let $n,k \in \N, n \geq k, p \in \pz$, all
fixed. Using $e_p (n!)=\sum_{i=1}^{\infty} [\frac{n}{p^i}]$ and
$[x+y]-[x]-[y] \in \{0,1\} \text{ } (\forall x,y \in \R)$ we
conclude:

$$ p | \binom{n}{k}$$
$$ \Leftrightarrow  \sum_{i=1}^{\infty}
\left(\left[\frac{n}{p^i}\right]-\left[\frac{k}{p^i}\right]-\left[\frac{n-k}{p^i}\right]\right) >0$$
$$\Leftrightarrow \exists i \in \N: \left[\frac{n}{p^i}\right]-\left[\frac{k}{p^i}\right] - \left[\frac{n-k}{p^i}\right] = 1 $$
\begin{equation}\Leftrightarrow \exists i \in \N, \exists j,f \in \N_0 :
\left[\frac{k}{p^i}\right] = j \text{ and }
\left[\frac{n-k}{p^i}\right]= f \text{ and }
\left[\frac{n}{p^i}\right]=f+j+1 \text .
\end{equation}

Now for $r \in \N_0, x \in \R$ and with the definition
$\frac{x}{0}=\infty$ we have $[\frac{x}{p^i}]=r \Leftrightarrow p^i
\in (\frac{x}{r+1}, \frac{x}{r}]$. Therefore (2) is equivalent to
the condition that $\exists i \in \N, \exists j,f \in \N_0 $, such
that with the following notation we have $p^i \in (a_1,a_2] \cap
(b_1,b_2] \cap (c_1,c_2] (=: A(j,f)=: A):$
\begin{equation}
p^i \in \left(\frac{k}{j+1}, \frac{k}{j}\right] =: (a_1,a_2],
\end{equation}
\begin{equation}
p^i \in \left(\frac{n}{f+j+2},\frac{n}{f+j+1}\right] =:(b_1,b_2]
\end{equation}
\begin{equation}
p^i \in \left(\frac{n-k}{f+1},\frac{n-k}{f}\right] =: (c_1, c_2].
\end{equation}
\\
We compute the intersection $A$ for fixed $j,f \in \N_0$: We assume
that $A \neq \emptyset $. Hence, we must have $a_1 < b_2$ and $c_1 <
b_2$. This yields (we omit the easy calculations):

$$
a_1 < b_2 \Leftrightarrow f < \left(\frac{n}{k} -
1\right)\left(j+1\right) \Leftrightarrow b_2 < c_2, $$
\begin{equation} c_1 < b_2 \Leftrightarrow f > \left(\frac{n}{k} -
1\right)j - 1 \Leftrightarrow b_2 < a_2 .
\end{equation}

Since $\max \{a_1,b_1,c_1\}<b_2=\min \{a_2,b_2,c_2\}$ we have $A
\neq \emptyset$ and the upper bound of the interval $A$ is $b_2$.
For the lower bound, we have to distinguish two cases.

Case 1: $f \leq (\frac{n}{k} - 1)(j+1) - 1 $ . We have $f \leq
(\frac{n}{k} - 1)(j+1) - 1 \Leftrightarrow a_1 \leq b_1
\Leftrightarrow b_1 \leq c_1$. So we get $A=(c_1,b_2]$.

Case 2: If $f > (\frac{n}{k} - 1)(j+1) - 1$, then using the other
bound from (6) we get $f = [(\frac{n}{k} - 1)(j+1)]$ and
$(\frac{n}{k} - 1)(j+1) \notin \Z$. The last condition is equivalent
to $n(j+1) \neq 0 \mod k$.
According to the calculations in case 1, in case 2 we get $c_1 < b_1$ and $b_1 < a_1$. So we get $A=(a_1,b_2]$.\\

Summarizing we calculated
$$A= \left\lbrace \begin{array}{ll} (c_1,b_2] \neq \emptyset&, j \in \N_0 \text{ and } f \in \left(\left(\frac{n}{k} -
1\right)j - 1, (\frac{n}{k} - 1)(j+1) - 1 \right] \cap \N_0 \\
 (a_1,b_2] \neq \emptyset &,j \in \N_0 \text{ and } f =[(\frac{n}{k} - 1)(j+1)] \text{ and } n(j+1) \neq 0 \mod k\\ \emptyset&,else\\
\end{array} \right. .$$
Thus, we have proven equation (1) (after having shifted the index $j$).\\

It remains to show that the intervals $A=A(j,f)$ are disjoint: Let
$(j_1,f_1) \neq (j_2,f_2)$. Case 1: $j_1 \neq j_2$. Then (3) shows
that $A(j_1,f_1) \cap A(j_2,f_2) = \emptyset$ (recall that
$A(j_l,f_l)$ is the intersection of the three intervals, given by
(3)-(5)). Case 2: $j_1 = j_2, f_1 \neq f_2$. Now (4) (also (5))
shows that the intersection is empty. \Endproof \vspace{1cm}

\textbf{Examples} $$\{p \in \pz; p | \binom{2000}{1000} \}$$
$$=\pz \cap \left( (1000,2000] \cup (500, 666] \cup (333,400] \cup (250,285] \cup \dots
\right),$$\\
$$\{p \in \pz; p | \binom{2000}{800} \}$$
$$=\pz \cap \left( (1200,2000] \cup (800, 1000] \cup (600,666] \cup (400,500]
\cup (300, 333] \cup (266,285] \cup \dots \right) .$$\\

 \begin{description}
 \item[Theorem 1] \emph{Let $m,n \in \N , n \geq m$. We have:}

\begin{equation}\omega \binom {nk}{mk} = \sum_{j=1}^{\infty} \left( \pi
\left(\frac{nk}{j}\right) - \pi \left(\frac{(n-m)k}{j} \right)- \pi
\left(\frac{mk}{j} \right) \right) + O(\sqrt{k}) \text{ \emph{for} }
k \longrightarrow \infty  .\end{equation}

 \end{description}
\vspace*{0.5cm} \Proof With Lemma 1 (take only the intervals for
$i=1$, the rest is a subset of $[1,\sqrt{nk}] \Rightarrow$ error of
$O(\sqrt{k})$) we get:
$$\omega \binom {nk}{mk} = \sum_{j=1}^{\infty}
\sum_{f=[\frac{n}{m}(j-1)]-j+1 }^{[\frac{n}{m}j]- j - 1} \left( \pi
\left(\frac{nk}{f+j}\right) - \pi \left(\frac{(n-m)k}{f+1} \right)
\right)$$
\begin{equation} + \sum_{j=1,\text{ } nj \neq 0 \text{ } mod\text{ } m}^{\infty}
\left( \pi \left(\frac{nk}{[\frac{n}{m}j]}\right) - \pi
\left(\frac{mk}{j}\right) \right) + O(\sqrt{k}) \text{ \emph{for} }
k \longrightarrow \infty  .
\end{equation}\\
Now (7) follows from (8) since (note that for each fixed $k$ we have
finite sums):
$$\sum_{j=1}^{\infty} \sum_{f=[\frac{n}{m}(j-1)]-j+1
}^{[\frac{n}{m}j]- j - 1} \pi \left(\frac{(n-m)k}{f+1}\right) =
\sum_{j=1}^{\infty} \pi \left(\frac{(n-m)k}{j}\right),$$
$$
\sum_{j=1,\text{ } nj \neq 0 \text{ } mod\text{ } m}^{\infty} \left(
\pi \left(\frac{nk}{[\frac{n}{m}j]}\right) - \pi
\left(\frac{mk}{j}\right) \right) =\sum_{j=1}^{\infty} \left( \pi
\left(\frac{nk}{[\frac{n}{m}j]}\right) - \pi
\left(\frac{mk}{j}\right) \right),
$$
$$
\sum_{j=1}^{\infty} \sum_{f=[\frac{n}{m}(j-1)]-j+1
}^{[\frac{n}{m}j]- j - 1} \pi \left(\frac{nk}{f+j}\right) +
\sum_{j=1}^{\infty} \pi \left(\frac{nk}{[\frac{n}{m}j]}\right) =
\sum_{j=1}^{\infty} \pi \left(\frac{nk}{j}\right) .
$$

\Endproof \vspace{1cm}

 Paul Erd\ös mentioned in [Erd79] that
 \begin{equation}
 \omega \binom{n}{k} = (1 + o(1))\frac{\log \binom {n}{k}}{\log n}
 \text{ for } k > n^{1-o(1)}.
 \end{equation}
 A proof for this fact can be easily obtained if one looks at Erd\ös' proof for a
 weaker statement given in [Erd73, p. 53]. Since
 Erd\ös didn't explicitly write down a proof for (9), we do it
 (we formulate the Theorem just for the case in which we are interested).\\

 Before we proceed with the proof, note the interesting fact that it
is comparatively difficult to show that $\omega ((2n)!)= \frac
{2n}{\log 2n} + o(\frac{n}{\log n})$ (Prime Number Theorem), but
much easier to show e.g. $\omega (\frac{(2n)!}{n!n!})= 2\log 2 \frac
{n}{\log n} + o(\frac{n}{\log n})$  (Theorem 2).

 \begin{description}
 \item[Theorem 2] (\emph{Erd\ös, [Erd79]}) \emph{Let $m,n \in \N , n \geq m$. We have:}
  \begin{equation}
\omega \binom {nk}{mk} = \log \frac {n^n}{m^m
(n-m)^{n-m}} \frac{k}{\log k} + o\left(\frac{k}{\log k}\right) \text{ for } k \longrightarrow \infty.\\
\end{equation}
 \end{description}
\vspace*{0.5cm} \Proof The proof looks as follows: $$\binom {nk}{mk}
= \prod_{p \mid \binom {nk}{mk} } p^{e_p(\binom {nk}{mk})}
 \approx \prod_{p \mid \binom {nk}{mk} } (nk) = (nk) ^{ \omega \binom {nk}{mk}} .
 $$ Since one can easily compute $\log \binom {nk}{mk}$, the theorem follows.\\

So let $m,n \in \N$ be fixed with $n>m$ (for $n=m$ the Theorem is
obviously correct). Using $\log k! =\int_1^k \log t dt+ O(\log k) =
k \log k - k + O(\log k)$ we get
$$ \log \binom {nk}{mk} = \log (nk)! - \log (mk)! - \log ((n-m)k)!
$$
$$ = nk \log n - mk \log m - (n-m)k \log (n-m)
+ O(\log k)
$$
$$
= k \log \frac {n^n}{m^m (n-m)^{n-m}} +O(\log k).
$$
\\
We now need the crucial fact, that ([Her68],[Sta69],[Sch69]):
$$\forall p \in \pz, n,k \in \N , n \geq k : p^{e_p (\binom{n}{k})}
\leq n \text{ } .$$ (Proof: $e_p (\binom
{n}{k})=e_p(n!)-e_p(k!)-e_p((n-k)!)=\sum_{i=1}^{[\frac{\log n}{\log
p}]} ([\frac{n}{p^i}]-[\frac{k}{p^i}]-[\frac{(n-k)}{p^i}]) \leq
\sum_{i=1}^{[\frac{\log n}{\log p}]} 1 \leq \frac{\log n}{\log p}
\Rightarrow p^{e_p (\binom{n}{k})} \leq p^{\frac{\log n}{\log p}} =
n $). We therefore get:
$$ \binom {nk}{mk} \leq (nk)^{\omega \binom {nk}{mk}} $$ $$\Rightarrow
\omega \binom {nk}{mk} \geq \frac{\log \binom{nk}{mk}}{\log nk} =
\frac{k}{\log k} \log \frac {n^n}{m^m (n-m)^{n-m}} +
o\left(\frac{k}{\log k}\right).
$$
Thus, we have the $\geq$ of (10). As for the $\leq$ we proceed as
follows:\\

Let $\varepsilon \in (0,1), k \in \N$ and $f(k,\varepsilon):= \#
\{p>(nk)^{1-\varepsilon}; p \in \pz, p \mid \binom{nk}{mk}\}$. We
have

$$ \binom {nk}{mk} > ((nk)^{1-\varepsilon})^{f(k,\varepsilon)} $$
$$\Rightarrow
f(k,\varepsilon) < \frac{\log \binom{nk}{mk}}{(1-\varepsilon)\log
nk} =\frac{k}{(1-\varepsilon) \log k} \log \frac {n^n}{m^m
(n-m)^{n-m}} + o\left(\frac{k}{\log k}\right).$$
\\
Therefore we get
$$ \omega \binom {nk}{mk} \leq f(k,\varepsilon) + k^{1-\varepsilon}
< \frac{k}{(1-\varepsilon) \log k} \log \frac {n^n}{m^m (n-m)^{n-m}}
+ o\left(\frac{k}{\log k}\right).$$

Since the last inequality holds for all $\varepsilon \in (0,1)$ we
get
$$ \limsup_{k \rightarrow \infty} \frac{\omega \binom {nk}{mk}}{\frac{k}{\log k}}
\leq \log \frac {n^n}{m^m (n-m)^{n-m}}$$ and the proof is complete.
\Endproof
 \vspace*{1cm}
 \textbf{Remark} Since $\pi (x) = \frac{x}{\log x} +
o(\frac{x}{\log
 x})$, Theorem 2 shows that the binomial coefficients $\binom {nk}{mk}$
have a \textquotedblleft positive proportion\textquotedblright  of
all prime numbers $\leq nk$ as divisors.\\
Note that if $k=o(n)$ ($k=k(t), n=n(t), k \leq n, \lim_{t
\rightarrow \infty} k(t)= \infty$), then $\omega \binom
{n}{k}=o(\frac{n}{\log
n}) \text{ for } t \rightarrow \infty$.\\
\Proof Like in the second part of Theorem 2's proof, we get for an
$\varepsilon \in (0,1)$:
$$\omega \binom{n}{k} <\frac{1}{1-\varepsilon} \frac{1}{\log n}
\log \frac{n^n}{k^k (n-k)^{n-k}} + o \left(\frac{n}{\log n}\right).
$$ Define $l=\frac{n}{k}$, then $l \longrightarrow \infty$ as $t
\rightarrow \infty$ and it follows, that $\frac{1}{n} \log
\frac{n^n}{k^k (n-k)^{n-k}}=\log \frac{l}{l-1} + \frac{1}{l} \log
(l-1) \longrightarrow 0$ whence we get the statement. \Endproof
\vspace{1cm}
 \textbf{Examples} Some concrete values which one yields
from Theorem 2 are (note that
$\lim \omega \binom{k}{\frac{2}{5} k} / \frac{k}{\log k}= \lim \omega \binom{5k}{2k} / \frac{5k}{\log k}$)\\
$\omega \binom{k}{k/2} / \frac{k}{\log k} \rightarrow 0.69\dots; $
$\omega \binom{k}{k/3} / \frac{k}{\log k} \rightarrow 0.63\dots; $
$\omega \binom{k}{k/4} / \frac{k}{\log k} \rightarrow 0.56\dots; $\\
$\omega \binom{k}{k/5} / \frac{k}{\log k} \rightarrow 0.50\dots; $
$\omega \binom{k}{k/10} / \frac{k}{\log k} \rightarrow 0.32\dots; $
$\omega \binom{k}{k/100} / \frac{k}{\log k} \rightarrow 0.05\dots; $\\
$\omega \binom{k}{k \frac{2}{5}} / \frac{k}{\log k} \rightarrow
0.67\dots$. \vspace{1cm}

By combining Theorem 1 and 2 we get for each pair $n,m \in \N, n>m$
an identity of the form $\alpha_{m,n} \frac{x}{\log x} +
o(\frac{x}{\log x}) = \sum_i (-1)^{a_{m,n}(i)}\pi
(\frac{x}{b_{m,n}(i)})$. We can sum up these identities for
different pairs $(n,m)$ and therefore get new identities. It turns
out that one can get these resulting identities (and even some more)
significantly easier for $\psi(x)$ and this has been mentioned quite
often, e.g. in [Land09, pp. 71-95]:
\begin{description}
 \item[Theorem 3] (\emph{e.g. [Land09, pp. 71-95]}) \emph{Let $l,j \in \N, m_1, \dots, m_l,n_1, \dots, n_j \in \N$
with $\sum_{i=1}^{l} m_i=\sum_{i=1}^{j} n_i$. We have:}
  \begin{equation}
\log \frac{(n_1k)! \dots (n_jk)!}{(m_1k)! \dots (m_lk)!} =
k \log \frac{n_1^{n_1} \dots n_j^{n_j}}{m_1^{m_1} \dots m_l^{m_l}} + O(\log k),\\
\end{equation}
\begin{equation}
\log \frac{(n_1k)! \dots (n_jk)!}{(m_1k)! \dots (m_lk)!} =
\sum_{i=1}^{\infty} \left( \psi \left(\frac{n_1k}{i} \right) + \dots + \left(\frac{n_jk}{i} \right)
- \psi \left(\frac{m_1k}{i} \right) - \dots - \left(\frac{m_lk}{i} \right) \right).\\
\end{equation}
 \end{description}
\vspace*{0.5cm} \Proof (11) follows from $\log k! = k \log k - k +
O(\log k)$.\\
As for (12): Put $ \Lambda (n) = \log p \text{ if } n=p^j \text{ }
(p \in \pz, j \in \N)$ and $\Lambda(n) = 0$ otherwise. Then
$$\log k! = \sum_{i=1}^k \log i = \sum_{i=1}^k \sum_{n | i} \Lambda
(n) = \sum_{n=1}^k \left[\frac{k}{n}\right] \Lambda (n) =
\sum_{n=1}^k \sum_{i \leq k/n} \Lambda (n) =$$
$$ \sum_{i=1}^k \sum_{n \leq k/i} \Lambda (n) = \sum_{i=1}^k \psi
\left(\frac{k}{i}\right)=
 \sum_{i=1}^{\infty} \psi \left(\frac{k}{i}\right) $$
 which proves (12). \Endproof
\vspace{1cm}
  \textbf{Remark} In Theorem 3 one can replace $\log (xk)!$ $(x\in \N)$ by
   $L(xk)=\sum_{i\leq xk} \log i$ $(x\in \R)$ and the
  stated identities remain correct.\\

\section{Applications}
In Corollary 1 we show how to use the identities obtained by
combining Theorem 1 and 2 in order to get some bounds for $\pi (x)$.
In the remark after the Corollary we note why a refinement of the
method used in Corollary 1 is not suitable for proving the Prime
Number Theorem.

Corollary 2 is a statement which includes the following well-known
result due to \emph{Chebychev}:
\begin{equation}\lim_{x\rightarrow \infty} \frac{\pi (x)}{x/\log x}
\text{ exists }
    \Leftrightarrow \lim_{x\rightarrow \infty} \frac{\pi (x)}{x/\log x} = 1.
\end{equation}

The paper ends with Corollary 3, which - using Theorem 1, 2 and the
Prime Number Theorem - allows one to give a series which converges
to $\log k$, thus generalizing the identity $\log 2 =
\sum_{j=1}^\infty \frac{(-1)^{j+1}}{j}$. Note that \textquotedblleft
elementary\textquotedblright proofs exist for this fact, which don't
make use of such \textquotedblleft heavy\textquotedblright Theorems
like the Prime Number Theorem (e.g. [KicGoe98],[Les01]).
Nevertheless, it might be interesting to see how this series occurs
again in connection
with the prime factors of certain binomial coefficients.\\

Corollary 1 gives some bounds for $\pi (x)$. Note that exactly the
same method and(!) numbers have been used often for doing the same
with $\psi (x)$ ([IwaKow04, p.33], [Land09, pp. 87-91]). In the
latter case Theorem 3 (with $n_1=30, n_2=1, m_1=15, m_2=10, m_3=6$)
is used instead of Theorem 1 and 2. We just included Corollary 1 and
its proof for the purpose of completeness.
\begin{description}
 \item[Corollary 1]
\begin{equation} 0.92 \frac{x}{\log x}+ o\left(\frac{x}{\log x}\right) < \pi (x) <
1.11 \frac{x}{\log x} + o\left(\frac{x}{\log x}\right).
\end{equation}
 \end{description}
\Proof With Theorem 1 (for $n=3,4,6$ and $m=1$; let $k \in 60 \N$)
we get
\\$\omega \binom{k/2}{k/6} = \pi(\frac{k}{2})
- \pi(\frac{k}{3}) + \pi(\frac{k}{4}) - \pi(\frac{k}{6}) + \pi(\frac{k}{8}) - \pi(\frac{k}{9}) + \pi (\frac{k}{10}) \mp ... + O(\sqrt{k})$\\
 $\omega \binom{k/3}{k/12} = \pi(\frac{k}{3}) -
\pi(\frac{k}{4}) + \pi(\frac{k}{6}) - \pi(\frac{k}{8}) + \pi(\frac{k}{9}) - \pi(\frac{k}{12}) \pm ... + O(\sqrt{k})$\\
$- \omega \binom{k/10}{k/60} =- \pi(\frac{k}{10}) +
\pi(\frac{k}{12}) - \pi(\frac{k}{20}) + \pi(\frac{k}{24}) - \pi(\frac{k}{30}) \pm ... + O(\sqrt{k})$\\
Sum up the last three identities to obtain \begin{equation}\omega
\binom{\frac{k}{2}}{\frac{k}{6}} + \omega
\binom{\frac{k}{3}}{\frac{k}{12}} - \omega
\binom{\frac{k}{10}}{\frac{k}{60}} = \pi\left(\frac{k}{2}\right) -
\pi\left(\frac{k}{12}\right) + \pi\left(\frac{k}{14}\right) -
\pi\left(\frac{k}{20}\right) \pm ... + O(\sqrt{k}) .\end{equation} A
further analysis shows that the sign on the right side of (3) is
alternating:\\
Let us associate to the sum on the right side of (3) the sequence
$(a_n)_{n\in \N}=(0,1,0,0,0,0,0,0,0,0,0,-1,0,1,\dots)$ where the
value $a_n$ shows how much times $\pi (\frac{x}{k})$ is added in the
sum. By Theorem 1 one sees that $(a_n)$ has period $60$ and we have:
$$a_n= \left\lbrace \begin{array}{ll} -1&,n \equiv \{12,20,24,30,36,40,48,60\} \mod 60 \\
 1&,n \equiv \{2,14,22,26,34,38,46,58\} \mod 60\\ 0&,else\\
\end{array} \right. .$$
\\
Thus, the sum on the right side of (3) is alternating. That is
important, since - because of the monotonicity of $\pi(x)$ - we can
then estimate

\begin{equation}\pi\left(\frac{k}{2}\right) - \pi\left(\frac{k}{12}\right) + O(\sqrt{k}) \leq \omega
\binom{\frac{k}{2}}{\frac{k}{6}} + \omega
\binom{\frac{k}{3}}{\frac{k}{12}} - \omega
\binom{\frac{k}{10}}{\frac{k}{60}}  \leq \pi\left(\frac{k}{2}\right)
+ O(\sqrt{k}).
\end{equation}

On the other hand Theorem 2 gives us
$$ \omega
\binom{\frac{k}{2}}{\frac{k}{6}} + \omega
\binom{\frac{k}{3}}{\frac{k}{12}} - \omega
\binom{\frac{k}{10}}{\frac{k}{60}} = \left(\frac{0.6365\dots}{2} +
\frac{0.5623\dots}{3} - \frac{0.4505\dots}{10} \right) \frac{k}{\log
k} + o\left(\frac{k}{\log k}\right)
$$
\begin{equation}
= 0.460 \dots \frac{k}{\log k} + o\left(\frac{k}{\log k}\right) .
\end{equation}

If we combine (4) and (5) and use the well-known bound
$\pi(\frac{k}{12}) \leq \frac{2}{12} \frac{k}{\log k}$, we get:
$$
0.92 \frac{x}{\log x}+ o\left(\frac{x}{\log x}\right) < \pi (x) <
1.26 \frac{x}{\log x} + o\left(\frac{x}{\log x}\right).
$$
If we use the last inequality for estimating $\pi(\frac{k}{12})$ and
do the last step again then we get the upper bound \textquotedblleft
$1.135$\textquotedblright. Do the last step once more to get the
upper bound \textquotedblleft $1.11$\textquotedblright. \Endproof
\vspace{1cm}

 \textbf{Remark 1} The statement of Corollary 1 can be proven also
 for $\psi(x)$ instead of $\pi(x)$ ([IwaKow04, p. 33], [Land09, pp. 87-91]).
 Just use Theorem 3 with $n_1=30, n_2=1, m_1=15, m_2=10,
 m_3=6$ to get the analogue of equation (3) for $\psi(x)$. Then
 continue as in the proof above.\\

 \textbf{Remark 2}
 Looking at the proof of Corollary 1, one immediately tries to add
 more $\omega \binom{k/l_i}{k/m_i}$'s in (3) in order to sharpen the
 estimates (hopefully to $1-\varepsilon < \frac{\pi(x) \log x}{x} < 1 + \varepsilon$ for
 each $\varepsilon >0$ $\Rightarrow$ Prime Number Theorem). In the
 case of (3) one could use $\omega \binom{k/12}{k/84}$ to get
 \begin{equation}\omega
\binom{\frac{k}{2}}{\frac{k}{6}} + \omega
\binom{\frac{k}{3}}{\frac{k}{12}} - \omega
\binom{\frac{k}{10}}{\frac{k}{60}} + \omega
\binom{\frac{k}{12}}{\frac{k}{84}} = \pi\left(\frac{k}{2}\right) -
 \pi\left(\frac{k}{20}\right) + \pi\left(\frac{k}{22}\right) + \pi\left(\frac{k}{26}\right) \mp \dots +
O(\sqrt{k}) .\end{equation}
 However, that is not as good as it may possibly look: Now the sum is not
 alternating anymore. While in Corollary 1 the \textquotedblleft error
 term\textquotedblright which appeared was $\pi(k/12)$, now the \textquotedblleft error
 term\textquotedblright is $\pi(k/20) + \pi(k/26)$. \\
\emph{Landau} explains in [Land09, pp. 597-598, 941] why an
extension of the above method (\emph{Landau} refers to the analogue
identities which one gets for $\psi(x)$, see Theorem 3) is not
suitable for proving the Prime Number Theorem (which we will call
from now on PNT). We will give a bit more details then he gives: Let
us shortly sketch how the method would have to be systematized if
one wants to deduce PNT. Let us therefore use Theorem 3 instead of
Theorem 1 and 2. \emph{Having in mind} what we did in Corollary 1,
we start with
\begin{equation}\log n!= \psi\left(\frac{n}{1}\right) + \psi\left(\frac{n}{2}\right) +
\psi\left(\frac{n}{3}\right) + \dots.\end{equation}

If - using the method of Corollary 1 - we want to get some bounds
$a,A$ with $a + o(1)< \psi(x) / x < A + o(1)$ and $A-a \leq
\varepsilon$, then it is necessary (though not sufficient) to
eliminate from the sum (7) the expressions $\psi(\frac{n}{2}),
\psi(\frac{n}{3}), \dots, \psi(\frac{n}{k_0})$ where $k_0 \in
\N,\frac{1}{k_0} \leq \varepsilon$. Therefore we must(!) subtract in
equation (7) the expressions $\log (\frac{n}{2})!, \log
(\frac{n}{3})!, \log (\frac{n}{5})!, \dots, \log
(\frac{n}{p_{k_0}})! \text{ } (p_i=\max \{p \in \pz; p \leq i\})$.
But then we have subtracted the expressions $\psi(\frac{x}{p_i
p_j})$ twice, thus we have to add these once. This has to be
continued and the outcome of the procedure is ($k_1 \geq k_0$; we
can assume $\frac{n}{p_{i1} \dots p_{ij}}, \frac{n}{l_0} \in \N$; we
also divide the equation through $n$)

\begin{equation} \frac{1}{n} \left( \psi (n) \pm
\psi\left(\frac{n}{k_1}\right) \pm \dots \right)=\end{equation} $$
\frac{1}{n} \left( \log \left(n! \cdot \prod_{p_{i1} p_{i2} \leq
k_0} \left(\frac{n}{p_{i1} p_{i2}}\right)! \cdot \dots \right) -
\log \left( \prod_{p_{i1}\leq k_0} \left(\frac{n}{p_{i1}}\right)!
\cdot \prod_{p_{i1} p_{i2} p_{i3} \leq k_0} \left(\frac{n}{p_{i1}
p_{i2} p_{i3}}\right)! \cdot \dots \right) \right).$$

Now for being able to work with (8) we have to add/subtract another
$\log (\frac{n}{l_0})!$ for an $l_0\in \Q$, so that the condition
$\sum_{i=1}^{l} m_i=\sum_{i=1}^{j} n_i$ of Theorem 3 is met. The
value of $l_0$ is determined by $k_0$:
$\frac{1}{l_0}=|\sum_{d=1}^{k_0} \frac{\mu(d)}{d}|$. Since we want
to get bounds $a,A$ with $A-a\leq \varepsilon$ it is necessary
(again not sufficient) that there exists a $k_0$ with
$\frac{1}{k_0}\leq \varepsilon$ such that
$\frac{1}{l_0}=|\sum_{d=1}^{k_0} \frac{\mu(d)}{d} |\leq
\varepsilon$. Since it is an \textquotedblleft
elementary\textquotedblright theorem that $\liminf \sum_{d=1}^{x}
\frac{\mu(d)}{d} \leq 0 \leq \limsup \sum_{d=1}^{x}
\frac{\mu(d)}{d}$ [Land09, pp. 583-584] we can get a suitable $k_0$
and $l_0$ (although new problems are likely to arise if $k_0$ is too
big). BUT there is a bigger problem: For our method to be
successful, the right side of (8) must tend to a limit. According to
Theorem 3 the right side of (8) is - apart from an error of $o(1)$ -
equal to

$$ \log  \left(\prod_{p_{i1} p_{i2} \leq k_0}
 \left(\frac{1}{p_{i1} p_{i2}}\right)^{\frac{1}{p_{i1} p_{i2}}} \cdot \dots \right) -
 \log  \left( \cdot \prod_{p_{i1} \leq k_0}
 \left(\frac{1}{p_{i1}}\right)^{\frac{1}{p_{i1}}} \cdot \prod_{p_{i1} p_{i2} p_{i3} \leq k_0}
 \left(\frac{1}{p_{i1} p_{i2} p_{i3}}\right)^{\frac{1}{p_{i1} p_{i2} p_{i3}}} \cdot \dots \right) =$$

 $$\sum_{d \leq k_0} \mu(d) \log \left( \frac{1}{d}^{\frac{1}{d}} \right)=- \sum_{d \leq k_0} \mu(d) \frac{\log d}{d}. $$

 Now, the statement that the latter sum converges for $k_0 \rightarrow
 \infty$ lies deeper then PNT as it is explained in [Land09, pp. 598-604, 941].
Let us summarize: For proving PNT via the above method one must
ensure (apart from other things) that $\sum_{d = 1}^{\infty} \mu(d)
\frac{\log d}{d}$ is convergent and once we know that, there is
already an elementary argument for
deducing PNT.\\
For more information on PNT see [BatDia96], [BatDia69], [Land09].\\
 \vspace{0.5cm}

 \begin{description}
 \item[Corollary 2] \emph{We have}
  \begin{equation} \forall k \in \N
\lim_{x\rightarrow \infty} \frac{\pi (kx)}{\pi (x)} \text{ exists }
\Leftrightarrow \lim_{x\rightarrow \infty} \frac{\pi (x)}{x/\log
x}=1 .\end{equation}
 \end{description}

\Proof
\\ \textquotedblleft $\Leftarrow$\textquotedblright $ \checkmark$ \\
     \textquotedblleft$\Rightarrow$\textquotedblright We will derive the statement from the following equation:
    \begin{equation} \log 2 \frac{x}{\log x} = \sum_{i=1}^{\infty} (-1)^{i+1} \pi \left(\frac{x}{i}\right) + o\left(\frac{x}{\log
x}\right) .\end{equation}
    Proof of (10): Using Theorem 1 and 2 for $n=2, m=1$, we
     derive the following identity:
$$\log 2 \frac{2k}{\log 2k} + o\left(\frac{k}{\log k}\right) = \pi (2k) - \pi \left(\frac{2k}{2}\right) + \pi
\left(\frac{2k}{3}\right) - \pi \left(\frac{2k}{4}\right) \pm \dots
 = \sum_{i=1}^{\infty}
    (-1)^{i+1} \pi \left(\frac{2k}{i}\right).$$\\
    Now, since for $|x-2k| \leq 2$: $\frac{2k}{\log 2k} - \frac{x}{\log x} = O(1)$ and $\sum_{i=1}^{\infty}
    (-1)^{i+1} \pi (\frac{2k}{i}) - \sum_{i=1}^{\infty}
    (-1)^{i+1} \pi (\frac{x}{i}) = \sum_{i=1}^{[\sqrt{2k}]}
    (-1)^{i+1} (\pi (\frac{2k}{i}) - \pi (\frac{x}{i}) ) + O(\sqrt{k}) = O (\sqrt{k}) $,
    we get what we wanted.
\\

    By assumption we can define $\forall k \in \N$:
    $\alpha_k :=\lim_{x \rightarrow \infty} \frac{\pi(kx)}{\pi(x)} \text{ }( \in \R ).$\\
    $  \sum_{i=1}^{\infty} (-1)^{i+1}
\frac{1}{\alpha_i} $ is convergent: $(\frac{1}{\alpha_i})_{i \in
\N}$ is monotonously decreasing. Furthermore $\alpha_2 > 1$,
otherwise $\alpha_2=1 \Rightarrow \alpha_{2^n} = 1 \text{ } \forall
n \in \N$, what would contradict $\frac{x}{\log x} \ll \pi(x) \ll
\frac{x}{\log x}$ (see Corollary 1). Altogether we get that
$(\frac{1}{\alpha_i})_{i \in \N}$ is a monotonously decreasing
zero-sequence ($\alpha_{2^i}=(\alpha_2 )^i \text{ } \forall i \in
\N)$, hence the mentioned series is convergent.

Now let $a_n := \sum_{i=1}^n (-1)^{i+1} \frac {1}{\alpha_i}$ and $a
:=\sum_{i=1}^{\infty} (-1)^{i+1} \frac {1}{\alpha_i}$. Take an
arbitrary $\varepsilon
> 0$ and choose $n_0 \in \N$, such that $\forall n \geq n_0 : |a-a_n| <
\varepsilon$. Using the monotonicity of $\pi (x)$ we get:
$$\frac{\sum_{i=1}^{\infty} (-1)^{i+1} \pi (\frac{x}{i}) }{\pi (x)}
\leq \sum_{i=1}^{2 n_0 + 1} (-1)^{i+1} \frac{\pi (\frac{x}{i}) }{\pi
(x)} \rightarrow \sum_{i=1}^{2 n_0 + 1} (-1)^{i+1}
\frac{1}{\alpha_i} \leq a + \varepsilon \text{ for } x \rightarrow
\infty \text{ as well as }$$
$$\frac{\sum_{i=1}^{\infty} (-1)^{i+1} \pi (\frac{x}{i}) }{\pi (x)}
\geq \sum_{i=1}^{2 n_0} (-1)^{i+1} \frac{\pi (\frac{x}{i}) }{\pi
(x)} \rightarrow \sum_{i=1}^{2 n_0} (-1)^{i+1} \frac{1}{\alpha_i}
\geq a - \varepsilon \text{ for } x \rightarrow \infty .$$

Thus we have $$\sum_{i=1}^{\infty} (-1)^{i+1} \pi
\left(\frac{x}{i}\right) = a \pi (x) + o\left(\frac{x}{\log
x}\right).$$

If we put the last equation into (10) we get  $$ \pi (x) =
\frac{\log 2}{a} \frac{x}{\log x} + o\left(\frac{x}{\log
x}\right).$$

We put the last identity into the definition of the $\alpha_k$ and
get
$$\forall k \in \N: \text{ } \alpha_k=k,$$
whence $$a=\sum_{i=1}^{\infty} (-1)^{i+1} \frac
{1}{\alpha_i}=\sum_{i=1}^{\infty} (-1)^{i+1} \frac {1}{i}=\log 2.$$
\Endproof \vspace{1cm}

 \textbf{Remark 1} Note that the condition
in (9, left side) is weaker then the one in (1, left side): (9, left
side) follows immediately from (1, left side), but the converse
doesn't follow instantly: E.g. let $f(x)$ be a function with
$f(x)=(x/\log x) (1 + \alpha \sin (\log \log x)) + o (x/\log x)
\text{ with } \alpha \in \R_{\neq 0} $.
    If we replace $\pi (x)$ by $f(x)$ then (9, left side) is true,
    but (1, left side) obviously not.\\

\textbf{Remark 2} Some superficial attempts to prove the Prime
Number Theorem by proving (9, left side) failed (e.g. by using an
explicit formula for $\pi(x)$ like the \textquotedblleft
Sieve of Eratosthenes\textquotedblright or  \textquotedblleft Meissels formula\textquotedblright).\\

 \textbf{Remark 3} Note that
Bertrands Postulate ($\pi (2n) - \pi(n)
> 0 $) can be verified immediately from equation (10) for $n \geq
n_0$, since $\pi(x) - \pi (x/2) + \pi(x/3) \geq \log 2 \frac{x}{\log
x} + o (\frac{x}{\log x})$, now use $\pi(x/3) <  \frac{2}{3}  \frac{x}{\log x}$.\\
On the other hand Bertrands Postulate for $n\geq n_0$ also follows
from the bounds in Corollary 1.\\

\textbf{Remark 4} The proof of Corollary 2 uses only the equation
(10) which can be obtained for $\psi(x)$ from Theorem 3 (take $n_1=2, m_1=1, m_2=1$).\\
\vspace{0.5cm}

There is an interesting generalization of the sum $\log 2 =
\sum_{j=1} \frac{(-1)^{j+1}}{j}$ for arbitrary $k \in \N$ instead of
$k=2$ (see for instance [KicGoe98], [Les01]).

With Theorem 1, 2 and the Prime Number Theorem we obtain a proof for
this generalization. Note that \textquotedblleft
elementary\textquotedblright proofs exist for Corollary 3, which
don't make use of such \textquotedblleft heavy\textquotedblright
Theorems like the Prime Number Theorem (e.g. [KicGoe98],[Les01]).
However, it might be interesting to see how these identities occur
again in connection with the prime factors of certain binomial
coefficients.

\begin{description}
 \item[Corollary 3] For all $k \in \N$ we have
\begin{equation}
\log k = \sum_{n=1}^{\infty} \left(\left(\sum_{i=1}^{k-1}
\frac{1}{nk-(k-i)}\right)-\frac{k-1}{nk}\right) \end{equation} $$= 1
+ \frac{1}{2}+ \dots \frac{1}{k-1} - \frac{k-1}{k}+ \frac{1}{k+1}+
\frac{1}{k+2} + \dots + \frac{1}{2k-1} - \frac{k-1}{2k} + \dots .
$$ \end{description}

\Proof
 Fix an arbitrary $n \in \N_{>1}$ and take $m=1$.
 Combine Theorem 1 and 2 and divide the derived identity through $\frac{k}{\log k}$.
 Then use the Prime Number Theorem to get

 (This doesn't fail because of any missing uniform convergence as it can
 be seen in the proof of Corollary 2, where - if we assumed the Prime Number Theorem as true -
 proved the following identity for the case $n=2$.):

\begin{equation} \log \frac{n^n}{(n-1)^{n-1}} = \sum_{j=0}^{\infty}
\sum_{i=1}^{n-1} \left( \frac{1}{j + i/n} - \frac{1}{j +i/(n-1)}
\right).
\end{equation}

Now having in mind that $\log \frac{k^k}{(k-1)^{k-1}} = k \log k -
(k-1) \log (k-1)$, it is an easy exercise in induction to obtain the
stated Corollary from (12). \Endproof \vspace{1cm}

\textbf{Examples} For instance we get
$$ \log 2 = \sum_{j=1}^{\infty} (-1)^{j+1} \frac{1}{j} ,$$
$$ \log 3 = \sum_{j=1}^{\infty} \frac{9j - 4}{(3j-2)(3j-1)3j}.$$
\vspace{1cm}
\section {References}

\small [BatDia69] Bateman P.T. and Diamond H.G., Asymptotic
distribution of Beurlings generalized prime numbers, Studies in
Number Theory Volume 6, W. J. LeVeque, Math. Assoc. America,
Washington, DC, 1969,
pp. 152-210\\

[BatDia96] Bateman P.T. and Diamond H.G., A Hundred Years of Prime
Numbers, The American Mathematical Monthly, Vol. 103, No.9, Nov.
1996, pp. 729-741\\

[Erd73] Erd\ös P., \Über die Anzahl der Primfaktoren von $\binom{n}{k}$, Arch. Math. 24, 1973, pp. 53-54 \\

[Erd79] Erd\ös P., Some unconventional problems in number theory,
Acta Math. Acad. Sc. Hung., t. 33, 1979, p. 71-80 (the theorem referred to in the text is on page 77 bottom)\\

[Fel91] Felgner U., Estimates for the sequence of primes, Elemente
der Mathematik, 1991, pp. 17 - 25
(see e.g. Theorem 4.2 and Lemma 5.1/5.2)\\

[Her68] Hering F., Eine Beziehung zwischen Binoialkoeffizienten
und Prim-zahlpotenzen, Arch. Math. 19, 1968, pp. 411 - 412 \\

[IwaKow04] Iwaniec H., Kowalski E., Analytic Number Theory, American Mathematical Society / Colloquium Publications Volume 53, 2004, p. 33 \\

[KicGoe98] Kicey C. and Goel S., A Series for ln k, American Math. Monthly, 105, 1998, pp. 552-554\\

[Land09] Landau E., Handbuch der Lehre von der Verteilung der
Primzahlen (with an Appendix by Bateman P.T.), Leipzig, Teubner, 1909\\

[Lan79] Langevin M., Facteurs premiers des coefficients binomiaux,
Seminaire Delange-Pisot-Poitou, Theorie des nombres,
tome 20, no. 2, 1978-1979, exp. no. 27, pp. 1-15 \\

[Les01] Lesko J., A Series for ln k, The College Mathematics
Journal, Vol. 32, No. 2., Mar. 2001, pp. 119-122\\

[Sch69] Scheid H., Die Anzahl der Primfaktoren in $\binom{n}{k}$,
Arch. Math. 20, 1969, pp. 581-582 \\

[Sta69] Stahl W., Bemerkung zu einer Arbeit von Hering, Arch. Math.
20, 1969, p. 580 \\

\end{document}